\documentclass[fleqn]{mat01}
\usepackage{times,mathtimy,amssymb,latexsym}
\begin{document}

\setcounter{page}{287} \firstpage{287}

\def\disp{\displaystyle}
\def\ZZ{\mathbb{Z}}
\def\CC{\mathbb{C}}
\def\calH{\mathcal{H}}
\def\lr{\left(}
\def\rr{\right)}
\def\lle{\left\{}
\def\rge{\right\}}
\def\varep{\varepsilon}
\def\cc#1{\{#1\}}
\def\pp#1{\|#1\|}
\def\ai{\overset{i}{\longrightarrow}}
\def\api{\overset{\pi}{\longrightarrow}}

\def\ao{\overset{0}{\longrightarrow}}
\def\axp{\overset{\times p^2}{\longrightarrow}}


\def\theor{\trivlist \item[\hskip \labelsep{\bf Theorem.}]}
\def\propo{\trivlist \item[\hskip \labelsep{\rm PROPOSITION.}]}

\title{On the cohomology of orbit space of free \pmb{${\ZZ}_{p}$}-actions\\ on lens spaces}

\markboth{Hemant Kumar Singh and Tej Bahadur Singh}{The orbit space of free
$\ZZ_p$-actions on lens spaces}

\author{HEMANT KUMAR SINGH and TEJ BAHADUR SINGH}

\address{Department of Mathematics, University of Delhi, Delhi~110~007, India\\
\noindent E-mail: tej$_{-}$b$_{-}$singh@yahoo.com}

 \volume{117}

\mon{August}

\parts{3}

\pubyear{2007}

\Date{MS received 11 November 2005; revised 12 August 2006}

\begin{abstract}
Let $G  = \ZZ_p$, $p$ an odd prime, act freely on a
finite-dimensional CW-complex $X$ with mod $p$ cohomology
isomorphic to that of a lens space $L^{2m-1} (p;q_1,\dots,q_m)$.
In this paper, we determine the mod $p$ cohomology ring of the
orbit space $X/G$, when  $p^2\nmid m$.
\end{abstract}

\keyword{Lens space; free action; cohomology algebra; spectral sequence.}

\maketitle

\section{Introduction}

Let $p$ be an odd prime and $m > 1$ an integer. Consider the
$(2m-1)$-sphere $S^{2m-1} \subset\CC \times \cdots\times \CC$
$(m$-times).  Given   integers $q_1,\dots,q_m$ relatively prime to
$p$, the map
$(\xi_1,\dots,\xi_m)\to(\zeta^{q_1}\xi_1,\dots,\zeta^{q_m}\xi_m)$,
where $\zeta=\hbox{e}^{2\pi\iota/p^2}$, defines a free action of
$G = \langle \zeta\rangle$ on $S^{2m-1}$. The orbit spaces of $G$
and the subgroup $N=\langle \zeta^p\rangle$ are the lens  spaces
$L^{2m-1}(p^2;q_1,\dots,q_m)$ and $L^{2m-1} (p;q_1,\dots,q_m)$,
respectively. Thus, we have a free action  of $\ZZ_p$ on $L^{2m-1}
(p;q_1,\dots,q_m)$ with the orbit space
$L^{2m-1}(p^2;q_1,\dots,q_m)$. By a mod $p$ cohomology lens space,
we mean a space $X$ whose \v{C}ech cohomology $H^*(X;\ZZ_p)$ is
isomorphic to that of a lens space $L^{2m-1}(p;q_1,\dots,q_m)$. We
will write $X\sim_p L^{2m-1}(p;q_1,\dots,q_m)$ to indicate this
fact. If $G=\ZZ_p$ acts on a mod $p$ cohomology lens space $X$,
then the fixed point set of $G$ on $X$ has been investigated by Su
\cite{4}. In this paper, we determine the cohomology ring (mod
$p$) of the orbit space $X/G$, when $G$ acts freely on $X$. The
following theorem is established.

\begin{theor}{\it
Let $G=\ZZ_p$ act freely on a finite-dimensional CW-complex $X\sim_p
L^{2m-1}(p;q_1,\dots,q_m)$. If $p^2\nmid m,$ then $H^*(X/G;\ZZ_p)$ is one of the
following graded commutative algebras$:$

\begin{enumerate}
\renewcommand\labelenumi{\rm (\roman{enumi})}
\leftskip .2pc
\item $\ZZ_p[x,y,z,u_1,u_3,\dots,u_{2p-3}]/I,$
where $I$ is the homogeneous ideal
\begin{align*}
&\hskip -1.25pc \langle x^2,y^p,z^n,u_hy-A_hxy^{(h+1)/2},u_hu_{2p-h},u_hu_{h'}-B_{hh'}xu_{h+h'-1}\\[.2pc]
&\hskip -1.25pc \quad\,-C_{hh'}y^{(h+h')/2},
u_hu_{h'}-B'_{hh'}zxu_{h+h'-2p-1}-C'_{hh'}zy^{(h+h'-2p)/2}\rangle \,,
\end{align*}
$m=np,$ $\deg x=1,$ $\deg y=2,$ $\deg z=2p,$ $\deg u_h =h,$ $\beta_p(x)=y,$ and $0 =
B_{hh'} = C_{hh'} = B'_{hh'} = C'_{hh'}$ when $h= h'$. $(\beta_p$ is the mod-$p$
Bockstein homomorphism associated with the sequence $0\to\ZZ_p\to\ZZ_{p^2}\to \ZZ_p\to
0)$.

\item $\ZZ_p [x,z]/\langle x^2,z^m),$
where $\deg x = 1$ and $\deg z = 2$.
\end{enumerate}}
\end{theor}

\section{Preliminaries}

In this section, we recall some known facts which will be used in the proof of our
theorem. Given a $G$-space $X$, there is an associated fibration $X\ai X_G\api B_G$,
and a map $\eta\hbox{\rm :}\ X_G\to X/G$, where $X_G=(E_G\times X)/G$ and $E_G\to B_G$
is the universal $G$-bundle. When  $G$ acts  freely on $X$, $\eta\hbox{\rm :}\ X_G\to
X/G$ is homotopy equivalance, so the cohomology rings $H^*(X_G)$ and $H^*(X/G)$ (with
coefficients in a field) are  isomorphic. To compute $H^*(X_G)$, we exploit  the
Leray--Serre spectral sequence of the map $\pi\hbox{\rm :}\ X_G\to B_G$. The
$E_2$-term of this spectral sequence is given by
\begin{equation*}
E_2^{k,l}\cong H^k (B_G;\calH^l(X))
\end{equation*}
(where $\calH^l(X)$ is a locally constant sheaf with stalk $H^l(X)$
and  group $G)$ and it converges to $H^*(X_G)$, as an algebra. The
cup product in $E_{r+1}$ is
induced from that in $E_r$ via the isomorphism $E_{r+1}\cong H^*(E_r)$.
When $\pi_1(B_G)$ operates trivially on $H^*(X)$, the system of
local coefficients is simple (constant) so that
\begin{equation*}
E_2^{k,l}\cong H^k (B_G)\otimes H^l(X)\,.
\end{equation*}
In this case, the restriction of the product structure in the spectral sequence to the
subalgebras $E_2^{*,0}$ and $E_2^{0,*}$ gives the cup products on $H^*(B_G)$ and
$H^*(X)$ respectively. The edge homomorphisms
\begin{align*}
&H^p(B_G)=E_2^{p,0}\to E_3^{p,0}\to \cdots\to E_{p+1}^{p,0}=E_\infty^{p,0}\subseteq
H^p (X_G) \ \ \mbox{and}\\[.2pc]
&H^q(X_G)\to E_\infty^{0,q}=E_{q+1}^{0,q}\subset \cdots\subset E_2^{0,q}=H^q(X)
\end{align*}
are the homomorphisms
\begin{equation*}
\pi^*\hbox{\rm :}\ H^p(B_G)\to H^p(X_G) \quad \hbox{and} \quad
\iota^* \hbox{\rm :}\ H^q(X_G)\to H^q(X),
\end{equation*}
respectively.

The above results about spectral  sequences can be found in \cite{3}. We also recall
that
\begin{equation*}
H^*(B_G;\ZZ_p) =\ZZ_p [s,t]/ (s^2) = \Lambda (s) \otimes \ZZ_p[t],
\end{equation*}
where $\deg s= 1$, $\deg t = 2$ and $\beta_p (s) = t$.

The following fact will be used without mentioning it explicitly.

\begin{propo}$\left.\right.$\vspace{.5pc}

\noindent {\it Suppose that $G = \ZZ_p$ acts on a finite-dimensional CW-complex space
$X$ with the fixed point set $F$. If $H^j(X;\ZZ_p) = 0$ for $j> n,$ then the inclusion
map $F\to X$ induces an isomorphism
\begin{equation*}
H^j(X_G;\ZZ_p)\to H^j(F_G;\ZZ_p)
\end{equation*}
for $j > n$ $($see Theorem~$1.5$ in Chapter~VII of \cite{1}).}
\end{propo}

\section{Proof}

The example of free action of $G=\ZZ_p$ on the lens space $L^{2m-1}(p;q_1,\dots,q_m)$
described in the introduction is a test case for the general theorem. All cohomology
groups in the proof should be considered to have coefficients in $\ZZ_p.$ Since
$\pi_1(B_G) = G$ acts trivially on $H^*(X)$, the fibration $X\ai X_G\api B_G$ has a
simple system of local coefficients on $B_G$. So the spectral sequence has
\begin{equation*}
E_2^{k,l}\cong H^k (B_G)\otimes H^l(X).
\end{equation*}
Let $a \in H^1(X)$ and $b\in H^2(X)$ be generators of the cohomology
ring $H^*(X)$.  As  there  are no fixed points  of  $G$  on  $X$,  the
spectral  sequence  does  not  collapse  at  the  $E_2$-term.
Consequently, we have either $d_2 (1\otimes a) = t \otimes 1$
and $d_2 (1 \otimes b) = 0$ or $d_2 (1\otimes a) = 0$ and $d_2 (1 \otimes b) =
t \otimes a$.

If  $d_2 (1 \otimes a) = 0$ and $d_2 (1 \otimes b) = t \otimes a$, then  we
have $d_2 (1 \otimes b^q) = q t \otimes ab^q$
and $d_2 (1\otimes ab^q)=0$
for $1\le q\le m-1$.
So $0 = d_2 [(1\otimes b^{m-1})\cup (1 \otimes b)]= m t\otimes ab^{m-1}$,
which is true iff $p\,|\,m$.
Suppose that $m = np$. Then
\begin{equation*}
d_2\hbox{\rm :}\ E_2^{k,l}\to E_{2}^{k+2,l-1}
\end{equation*}
is  an isomorphism if $l$ is even and $2p\nmid l$, and $d_2=0$ if
$l$ is odd or $2p\,|\, l$. So $E_3^{k,l} =E_2^{k,l}$ for all $k$
if $l=2qp$ or $2(q+1)p-1$, where $0\le q<n$; $k  =  0,1$  if  $l$
is odd and $2 p\nmid (l+1)$, and  $E_3^{k,l}=0$, otherwise. It is
easily seen that $d_3 =0$, for example, if $u\in
E_3^{0,2(q+1)p+1}$ and $d_3(u)=A[st\otimes ab^{(q+1)p-1}]
(A\in\ZZ_p)$, then, for $v=[t\otimes 1]\in E_3^{2,0}$, we have
$0=d_3 (u\cup v) = A [st^2\otimes ab^{(q+1)p-1}]\Rightarrow A =
0$. A similar argument shows that the differentials
$d_4,\dots,d_{2p-1}$ are all trivial. If
\begin{equation*}
d_{2p}\hbox{\rm :}\  E^{0,2p-1}_{2p}\to E_{2p}^{2p,0}
\end{equation*}
is also trivial, then
\begin{equation*}
d_{2p}\hbox{\rm :}\  E_{2p}^{k,l}\to E^{k+2p,l-2p+1}_{2p}
\end{equation*}
is  trivial  for  every $k$ and $l$, because  every  element  of
$E_{2p}^{k,2(q+1)p-1}$ can be written as the
product of an element of $E_{2p}^{k,2qp}$ by
$1\otimes ab^{p-1} \in E_{2p}^{0,2p-1}$ and
\begin{equation*}
d_{2p}\hbox{\rm :}\  E_{2p}^{k,2(q+1)p}\to E_{2p}^{k+2p,2qp+1}
\end{equation*}
is  obviously trivial. If $n = 1$, then $E_\infty = E_3$, where  the top and bottom
lines survive. This contradicts our hypothesis; so $n>1$. If $d_{2p+1} [1\otimes b^p]
= [st^p\otimes 1]$, then it can be easily verified that
\begin{align*}
&d_{2p+1} [1\otimes b^{qp}] = q[st^p\otimes b^{(q-1)p}]\quad \text{and}\\[.2pc]
&d_{2p+1}[1\otimes ab^{(q+1)p-1}] = q[st^p\otimes ab^{qp-1}]
\end{align*}
for $1\le q < n$. Consequently,
\begin{equation*}
0= d_{2p+1} [(1\otimes ab^{np-1})\cup (1\otimes b^p)] = n (st^p\otimes ab^{np-1}),
\end{equation*}
which is not true for $(n, p) = 1$. On the other hand, if
\begin{equation*}
d_{2p+1}\hbox{\rm :}\ E_{2p}^{0,2p}\to E_{2p}^{2p+1,0}
\end{equation*}
is trivial, then
\begin{equation*}
d_{2p+1}\hbox{\rm :}\ E_{2p}^{k,l} \to E^{k+2p+1,l-2p}_{2p}
\end{equation*}
is also trivial  for every $k$ and $l$, as above. Now $d_r=0$ for every $r > 2 p + 1$,
so several lines of the spectral sequence survive to infinity. This contradicts our
hypothesis. Therefore,
\begin{equation*}
d_{2p}\hbox{\rm :}\  E_{2p}^{0,2p-1}\to E_{2p}^{2p,0}
\end{equation*}
must be non trivial. Assume that $d_{2p}[1\otimes ab^{p-1}] =[t^p\otimes 1]$. Then
\begin{equation*}
d_{2p}\hbox{\rm :}\ E^{k,l+2p-1}_{2p}\to E^{k+2p,l}_{2p}
\end{equation*}
is an isomorphism for $l=2qp$, $0\le q < n$, and is trivial homomorphism for other
values of $l$. Accordingly, we have $E_\infty = E_{2p+1}$,  and hence
\begin{equation*}
H^j(X_G)= \begin{cases}
\ZZ_p, &j= 2qp, \ 2(q+1)p-1, \ 0\le q<n;\\
\ZZ_p\oplus\ZZ_p, & 2qp<j<2(q+1)p-1, \ 0 \le q <n; \ \text{and}\\
0, &j>2np-1.
\end{cases}
\end{equation*}
The elements $1 \otimes b^p\in E_2^{0,2p}$ and $1\otimes ab^{(h-1)/2} \in E_2^{0,h}$,
for $h = 1,3,\dots,2p-3$ are permanent cocyles; so they determine elements $z\in
E_\infty^{0,2p}$ and $w_h\in E_\infty^{0,h}$, respectively. Obviously, $\iota^* (z) =
b^p$, $z^n = 0$ and $w_hw_{h'}=0$. Let $x =\pi^* (s)\in E_\infty ^{1,0}$ and $y= \pi^*
(t) \in E_\infty ^{2,0}$. Then $x^2 = 0$, $y^p = 0$, and,  by  the naturality  of the
Bockstein homomorphism $\beta_p$, we have $\beta_p(x) = y$ and $yw_h=0$ but $xw_h\neq
0$. It follows that the total complex Tot $E_\infty^{*,*}$ is the graded commutative
algebra
\begin{align*}
\text{Tot }\,  E_\infty^{*,*} = \ZZ_p [x,y,z, w_1,w_3,\dots, w_{2p-3}]/\langle
x^2,y^p,z^n,w_hw_{h'},w_hy\rangle,
\end{align*}
where $h$, $h'=1,3,\dots,2p-3$. Choose $u_h\in H^h(X_G)$ representing $w_h$ for $h =
1,3,\dots,2p-3$. Then $\iota^*(u_h) = ab^{(h-1)/2}$, $u_h^2=0$ and $u_hu_{2p-h} = 0$.
It follows that
\begin{equation*}
H^*(X_G) = \ZZ_p [x,y,z,u_1,u_3,\dots,u_{2p-3}]/I,
\end{equation*}
where $I$ is the ideal generated by the homogenous elements
\begin{equation*}
\hskip -4pc x^2,y^p,z^n,yu_h-A_hxy^{(h+1)/2},u_h u_{2p-h},
u_hu_{h'}-B_{hh'}xu_{h+h'-1}-C_{hh'}y^{(h+h')/2}
\end{equation*}
and $u_hu_{h'}-B'_{hh'}zxu_{h+h'-2p-1}-C'_{hh'}zy^{(h+h'-2p)/2}$.

Here  $\deg  x = 1$, $\deg y = 2$, $\deg z = 2p$, $\deg u_h=h$  and,
when $h = h'$,
$0=B_{hh'}=C_{hh'}=B'_{hh'}=C'_{hh'}$.

If  $p\nmid m$, then we must have $d_2(1\otimes a) = t\otimes 1$, $d_2 (1 \otimes b)= 0$.
It can be easily observed that
\begin{equation*}
d_2\hbox{\rm :}\  E_2^{k,l}\to E_{2}^{k+2,l-1}
\end{equation*}
is a trivial homomorphism for $l$ even and an isomorphism for  $l$ odd. It is now easy
to see that $d_r = 0$ for every $r > 2$. So $E_\infty^{k,l} = E_3^{k,l} = \ZZ_p$ for
$k < 2$ and $l= 0,2,4,\dots,2m-2$. Therefore, we have
\begin{equation*}
H^j (X_G) =
\begin{cases}
\ZZ_p\,,& 0\le j\le 2m-1;\\
0\,,&\text{otherwise.}
\end{cases}
\end{equation*}
If  $x\in H^1(X_G)$ is determined by $s \otimes 1\in E_2^{1,0}$, then $x^2\in E_\infty^{2,0}=0$.
The multiplication by $x$
\begin{equation*}
x\cup (\cdot)\hbox{\rm :}\ E_\infty^{0,l}\to E_\infty^{1,l}
\end{equation*}
is an isomorphism for $l$ even. The element $1\otimes b \in E_2^{0,2}$ is a permanent
cocycle and determines an element $z \in E_\infty^{0,2} = H^2 (X_G)$. We have $\iota^*
(z) = b$ and $z^m = 0$. Therefore, the total complex Tot $E_\infty^{*,*}$ is the
graded commutative algebra
\begin{equation*}
\text{Tot}\,E_\infty^{*,*}=\ZZ_p[x,z]/\langle x^2,z^m\rangle\,.
\end{equation*}
Notice  that $H^j (X_G)$ is $E_\infty^{0,j}$ for $j$ even and $E_\infty^{1,j-1}$ for $j$  odd.
Hence,
\begin{equation*}
H^* (X_G) =\ZZ_p [x,z]/\langle x^2, z^m\rangle,
\end{equation*}
where $\deg x = 1$ and $\deg z = 2$. This completes the proof.
\hfill $\Box$

\section{Example}

We realize here the second case of our theorem. Recall that $G=\ZZ_p$ acts freely on
$L^{2m-1}(p;q_1,\dots,q_m)$  with the orbit space $L^{2m-1}(p^2;q_1,\dots,q_m)=K$. We
claim that
\begin{equation*}
H^*(K;\ZZ_p)=\ZZ_p[x,z]/\langle x^2,z^m\rangle,
\end{equation*}
where $\deg x=1, \ \deg z=2$. It is known that $K$ is a CW-complex with 1-cell of each
dimension $i=0,1,\dots,2m-1$ and the cellular chain complex of $K$ is
\begin{equation*}
0\to C_{2m-1}\ao C_{2m-2}\axp
C_{2m-3}\longrightarrow\cdots\longrightarrow C_2\axp C_1\ao C_0,
\end{equation*}
where each $C_i=\ZZ$. Accordingly, the co-chain complex of $K$  with coefficients in $\ZZ_p$  is
\begin{equation*}
0\to\ZZ_p\to\ZZ_p\to\cdots\to\ZZ_p\to \ZZ_p\to 0\,,
\end{equation*}
where each coboundary operator  is the trivial homomorphism. Therefore
\begin{equation*}
H^j(K;\ZZ_p)=
\begin{cases}
\ZZ_p\,,&\text{for $0\le j\le 2m-1$};\\
0\,,&\text{for $j\ge 2m.$}
\end{cases}
\end{equation*}
To determine the cup product in  $H^j(K;\ZZ_p)$, we first observe
that the inclusion $K^{(2i-1)}\to K^{(2i+1)}$ induces isomorphism
$H^j(K^{(2i-1)};\ZZ_p)\cong H^j(K^{(2i+1)};\ZZ_p)$ for  $j\le
2i-1$ so that we can identify them. For $j<2i-1$, this follows
from the cohomology exact sequence of the pair
$(K^{(2i+1)},K^{(2i-1)})$. The exact cohomology sequence of the
pairs $(K^{(2i+1)},K^{(2i)})$ and $(K^{(2i)},K^{(2i-1)})$ show
that $H^{2i-1}(K^{(2i+1)};\ZZ_p)\cong H^{2i-1}(K^{(2i)};\ZZ_p)$
and $H^{2i-1}(K^{(2i)};\ZZ_p)\cong H^{2i-1}(K^{(2i-1)};\ZZ_p)$;
the latter because the homomorphism
$H^{2i}(K^{(2i)},K^{(2i-1)};\ZZ_p)\to H^{2i}(K^{(2i)};\ZZ_p)$ is
surjective.

Now, we choose generators $x\in H^1(K;\ZZ_p)$ and $z\in H^2(K;\ZZ_p)$. Then obviously
$x^2=0$ and $z^m=0$. We can assume, by induction, that $z^i$  and $xz^i$ generate
$H^{2i}(K;\ZZ_p)$  and  $H^{2i+1}(K;\ZZ_p)$, respectively, for $i\le m-2$. Then, there
is an element $kxz^{m-2}$ such that $z\cup kxz^{m-2}=kxz^{m-1}$  generates
$H^{2m-1}(K;\ZZ_p)$ (see Corollary 3.39 of \cite{2}). We must have  $(k,p)=1$,
otherwise the order of $kxz^{m-1}$ would be less than $p$. Thus $xz^{m-1}$ generates
$H^{2m-1}(K;\ZZ_p)$, and this is true only if $z^{m-1}$ generates $H^{2m-2}(K;\ZZ_p)$.
Hence our claim.

\section{Remarks}

\begin{enumerate}
\renewcommand\labelenumi{\rm (\roman{enumi})}
\leftskip .2pc
\item It is clear from the proof of the theorem that if $p\nmid m$, then only
the second possibility of the theorem holds. Furthermore, if $X\sim_p
L^{2m-1}(p;q_1,\dots,q_m)$ and $\pi_1(X)=\ZZ_p$, then there exists a simply connected
space $Y$ with a free action of $\Delta=\ZZ_p$ such that $Y\sim_p S^{2m-1}$ and
$Y/\Delta\approx X$ (Theorems~3.11 and 2.6 of \cite{4}). If $G=\ZZ_p$ acts freely on
$X$, then the liftings of transformations (on $X$) induced by the elements of $G$ form
a group $\Gamma$ of order $p^2$ which acts freely on $Y$ and hence $\Gamma$ must be
cyclic. It is clear that $\Gamma$ contains the group $\Delta$ of deck transformations
of the covering $Y\rightarrow X$, and $G=\Gamma/\Delta$. So $X/G\approx Y/\Gamma$.
Since $Y\sim_p S^{2m-1}$, the mod $p$ cohomology algebra of $Y/\Gamma$ is a truncation
of $H^*(B_{\Gamma};\ZZ_p)$. Thus, in this case also, only the second possibility of
the theorem holds irrespective of the condition whether or not $p|m$.

\item We recall that a paracompact Hausdorff space $X$ is called finitistic if every
open covering of $X$ has a finite dimensional open refinement (see p.~133 of
\cite{1}). Our theorem and its proof go through for finitistic spaces.
\end{enumerate}

\section*{Acknowledgement}

We would like to thank the referee for his/her valuable suggestions which have
improved  exposition of the paper.

\end{document}